\documentclass[11pt,fleqn,leqno]{article}
\usepackage{amsmath}
\usepackage{amssymb}
\usepackage{amsfonts}
\usepackage{amsthm}
\usepackage{mathtools}
\mathtoolsset{%
above-intertext-sep = -1ex
}

\usepackage[utf8]{inputenc}
\usepackage[T1]{fontenc}

\usepackage[sans]{dsfont}

\usepackage[%
bb=fourier,
scr=euler,
frak=euler
]
{mathalfa}

\usepackage[font=small,labelfont=bf]{caption}
\usepackage{graphicx}
\graphicspath{{Figures/}} 
\usepackage{acronym}
\usepackage{latexsym}
\usepackage{paralist}
\usepackage{wasysym}
\usepackage{xspace}
\usepackage{framed}
\usepackage{palatino,pxfonts}
\usepackage{authblk}
\usepackage[numbers,sort&compress]{natbib}

\usepackage[dvipsnames,svgnames]{xcolor}
\colorlet{MyBlue}{DodgerBlue!75!Black}
\colorlet{MyGreen}{DarkGreen!95!Black}

\usepackage{hyperref}
\hypersetup{
colorlinks=true,
linktocpage=true,
pdfstartview=FitH,
breaklinks=true,
pdfpagemode=UseNone,
pageanchor=true,
pdfpagemode=UseOutlines,
plainpages=false,
bookmarksnumbered,
bookmarksopen=false,
bookmarksopenlevel=1,
hypertexnames=true,
pdfhighlight=/O,
urlcolor=MyBlue!60!black,linkcolor=MyBlue!70!black,citecolor=DarkGreen!70!black, 
pdftitle={},
pdfauthor={},
pdfsubject={},
pdfkeywords={},
pdfcreator={pdfLaTeX},
pdfproducer={LaTeX with hyperref}
}

\numberwithin{equation}{section}  
\usepackage[sort&compress,capitalize,nameinlink]{cleveref}
\crefname{example}{Ex.}{Exs.}

\crefrangeformat{equation}{\upshape(#3#1#4)\textendash(#5#2#6)}

\usepackage{autonum}

\newcommand{\dd}{\:d}

\newcommand{\eps}{\varepsilon}

\newcommand{\dif}{\dd}

\DeclareMathOperator*{\argmin}{argmin}

\renewcommand{\emptyset}{\varnothing}

\newcommand{\wlim}{\rightharpoonup}

\newcommand{\scrH}{\mathcal{H}}

\newcommand{\scrP}{\mathcal{P}}

\newcommand{\scrW}{\mathcal{W}}
\newcommand{\scrX}{\mathcal{X}}

\newcommand{\setH}{\mathsf{H}}

\newcommand{\Leb}{{\mathsf{Leb}}}

\newcommand{\R}{\mathbb{R}}

\newcommand{\N}{\mathbb{N}}

\DeclareMathOperator{\NC}{\mathsf{NC}}

\DeclareMathOperator{\VI}{VI}

\DeclareMathOperator*{\essinf}{ess\,inf}

\usepackage{algorithm2e}
\usepackage{algpseudocode}								

\theoremstyle{plain}
\newtheorem{theorem}{Theorem}

\newtheorem*{corollary*}{Corollary}
\newtheorem{lemma}[theorem]{Lemma}

\theoremstyle{definition}
\newtheorem{definition}[theorem]{Definition}
\newtheorem*{definition*}{Definition}
\newtheorem{assumption}{Assumption}

\renewcommand\qed{\hfill\small$\blacksquare$}

\theoremstyle{remark}

\newtheorem*{remark*}{Remark}
\newtheorem*{notation*}{Notational remark}

\numberwithin{theorem}{section}
\numberwithin{remark}{section}
\numberwithin{example}{section}

\DeclarePairedDelimiter{\abs}{\lvert}{\rvert}

\DeclarePairedDelimiter{\inner}{\langle}{\rangle}
\DeclarePairedDelimiter{\norm}{\lVert}{\rVert}

\title{Strong Convergence of Forward-Backward-Forward Methods for Pseudo-monotone Variational Inequalities with Applications to Dynamic User Equilibrium in Traffic Networks}
\date{\today}

\author[1]{\small Benoit Duvocelle}
\author[2]{\small Dennis Meier}
\author[1]{\small Mathias Staudigl}
\author[2]{\small Phan Tu Vuong}

\affil[1]{\footnotesize Maastricht University, Department of Quantitative Economics, P.O. Box 616, NL\textendash 6200 MD Maastricht, The Netherlands}
\affil[2]{\footnotesize Faculty of Mathematics, University of Vienna, Oskar-Morgenstern-Platz 1, A-1090 Vienna, Austria.}

\begin{document}

\maketitle

\begin{abstract}
In infinite-dimensional Hilbert spaces we device a class of strongly convergent primal-dual schemes for solving variational inequalities defined by a Lipschitz continuous and pseudomonote map. Our novel numerical scheme is based on Tseng's forward-backward-forward scheme, which is known to display weak convergence, unless very strong global monotonicity assumptions are made on the involved operators. We provide a simple augmentation of this algorithm which is computationally cheap and still guarantees strong convergence to a minimal norm solution of the underlying problem. We provide an adaptive extension of the algorithm, freeing us from requiring knowledge of the global Lipschitz constant. We test the performance of the algorithm in the computationally challenging task to find dynamic user equilibria in traffic networks and verify that our scheme is at least competitive to state-of-the-art solvers, and in some case even improve upon them.
\end{abstract}


\section{Introduction}
\label{sec:Intr}
Variational inequalities (VIs) are a flexible mathematical formulation of many equilibrium problems in engineering, machine learning, operations research and economics (see \cite{FacPan03} for a masterful survey of theory and applications of finite-dimensional VIs). Formulated on an infinite-dimensional real Hilbert space, variational inequalities also play a key role in the field of PDEs and optimal control, with important applications in imaging, differential equations, and network flows \cite{KinSta80,Hin08}. This paper is concerned with two major issues in studies of Hilbert-space valued variational inequalities. Our first aim to develop solution schemes with \emph{cheap iterations.} Second, we insist on numerical schemes guaranteeing \emph{strong convergence} of the generated trajectories. Both desiderata are motivated, among others, by developing provably convergent numerical solution techniques for a challenging problem in transportation science, the computation of \emph{dynamic user equilibria} in traffic networks.

\subsection{Dynamic user equilibrium}
Dynamic user equilibrium (DUE) is the most widely studied form of \emph{dynamic traffic assignment} (DTA), in which road travelers engage in a non-cooperative Nash game with departure time and route choices.  One characteristic feature of DTA is that it provides a "general equilibrium" model whose aim is to predict departure rates, departure times and route choices of travelers over a given time horizon. Exact DTA models are built on two layers: (i) a game-theoretic formulation of trip assignment, such as the dynamic extension of Wardrop's first principle \cite{War52}; (ii) a network flow model, which captures the physical relationships between entry and exit flows, junction flows, link delay and path delay. The latter is referred in the literature as \emph{dynamic network loading} (DNL). The DNL procedure is a manifestation of the physical principles of traffic flows, and various formulation of DNL exist in the literature, ranging from fluid models to differential equations; We refer the reader to the survey \cite{BreEMS14} and the book \cite{GarPic06} for an in-depth treatment of this important subject. We focus in this paper on the computation of DUE, leaving the network loading in the back. Section \ref{sec:DUE} gives a precise explanation how this division between the two levels works. A key challenge in the algorithmic approach to DUE is the usual lack of a closed-form expression of the \emph{delay operator}. The delay operator is the quantity of interest in DUE, since it informs us about the latencies on the individual paths of the traffic network. Indeed, as shown already in the seminal work \cite{Frie93}, the delay operator is the defining map in the VI approach of dynamic user equilibrium. However, without detailed information on this map, it is impossible to make a-priori monotonicity statements, which are crucial in the choice of numerical algorithms to solve the variational inequality. In fact, even if an explicit expression for the delay operator is available, it has been shown in \cite{MouSmi07} that strong monotonicity cannot hold for general networks and DNL models. Hence, any numerical algorithm guaranteeing \emph{strong convergence} under a-priori \emph{weak monotonicity} assumptions marks a breakthrough in the applicability of DUE as a predictive tool for traffic engineers.\\
The literature on DUE is huge, and naturally it is impossible to give a fair representation of all available results. We therefore only give a summary of those contributions which are the most related to our work. In finite-dimensions, the connection between VIs and traffic user equilibrium is classical (see e.g. \cite{FacPan03}). Once the user equilibrium problem is put into a dynamic setting, the natural model domain is the space of path-flows, which are assumed to be square-integrable functions satisfying a natural conservation condition. To our knowledge, the VI formulation of dynamic flows over time has been first presented in \cite{Frie93}. Departing from that work, the field has grown substantially, and various numerical schemes have been constructed to solve the resulting VI under different global regularity and Lipschitz continuity assumptions on the involved operator.\footnote{In terms of numerical analysis, these papers can thus be seen to follow the classical philosophy to \emph{first optimize, then discretize}.} A gradient projection method is studied in \cite{FriHanLiuSze13}. Weak convergence of this method is known if the operator is Lipschitz continuous and strongly monotone \cite{BauCom16}. As noted in \cite{HeLia02}, relaxing strong monotonicity assumptions could even lead to divergence of the algorithm. \cite{LoSze02} develop an alternating direction method under the assumption that the delay operator is \emph{cocoercive}. Sufficient for cocoercivity is Lipschitz continuity and monotonicity, so again we need to make rather restrictive global monotonicity assumption. Assuming weaker monotonicity conditions, the well-known extragradient scheme, due to Korpelevich and Antipin \cite{Kor76,Anti78}, has been employed in \cite{Long2013} to solve for DUE. In \cite{Vuong18} the weak convergence of the extragradient method is studied in some detail. A further drawback of the extragradient method is that it requires two costly projection steps at each iteration, making it a relatively unattractive method given our desire to have schemes with computationally cheap iterations. \cite{FriHanLiuSze13} also discuss a proximal-point algorithm, first studied by Martinet \cite{Mar70} and Rockafellar \cite{Roc76}, and the self-adaptive projection scheme of \cite{HanLo02}. Again, without assuming strong monotonicity, proximal-point methods are known to converge only weakly \cite{Guel91}, and the self-adaptive projection scheme has been introduced in \cite{HanLo02} in a finite-dimensional setting, making the distinction between weak and strong convergence meaningless. In light of the above survey, the following research question emerges:\\
\emph{
Can we develop a numerical algorithm with computationally cheap iterations and exhibiting strong convergence of the iterates under mild monotonicity assumptions?
}
\\
In this paper we provide an affirmative answer to this pressing question.

\subsection{Methodological contribution} 
Beside providing a new solution technique for a challenging infinite-dimensional equilibrium problem, we belief that our algorithmic design is also interesting from the point of view of numerical analysis. Specifically, our main result is achieved by constructing a novel iterative scheme for solving variational inequalities on a real Hilbert spaces by forcing strong convergence within the general framework of forward-backward-forward algorithms \cite{Tse00}. In infinite-dimensional settings strongly convergent iterative schemes are much more desirable than weakly convergent ones since it translates the physically tangible property that the energy $\norm{x^{k}-x^{\ast}}^{2}$ of the error between the iterate $x^k$ and a solution $x^{\ast}$ eventually becomes arbitrarily small. Of course, any numerical solution technique designed for solving a problem in infinite dimensions must be applied to a finite-dimensional approximation of the problem. Exactly in such situations strongly convergent methods are extremely powerful, because they guarantee stability with respect to discretizations. In fact, \cite{Guel91} demonstrated that strongly convergent schemes might even exhibit faster convergence rates as compared to their weakly convergent counterparts. It seems therefore fair to say that strong convergence is an extremely desirable property of solution schemes, with clearly observable physical consequences on the performance and stability of algorithms.\\
Our approach is closely related to \cite{GibTho18}, who develop a similar forward-backward-forward scheme in the setting of maximal monotone operators, and prove strong convergence of the iterates. We prove strong convergence under the weaker setting of pseudo-monotone VIs. The relaxation in the monotonicity of the operator is particularly valuable from the point of view of the motivating application for this paper, and we have discussed this issue extensively in the previous section. Very recently, the paper \cite{ThoChoVin19} developed a strongly convergent forward-backward-forward scheme, using hyperplane projections  \`{a} la Haugazeau \cite{Hau68} (see also \cite{SolSvai00,BauCom01} for an early application of this idea). Our numerical scheme is arguably much simpler, since we do not introduce any additional projection subroutines, but only simple extrapolations. \\
The rest of this paper is organized as follows. In Section \ref{sec:prelims} we introduce standard notation and concepts from variational analysis. Section \ref{sec:VI} describes the numerical scheme which we prove to be strongly convergent under Lipschitz continuity and weak-monotonicity assumptions only. An extension of the basic scheme to an adaptive algorithm is also discussed in that Section, showing that we can even get rid of Lipschitz continuity assumptions when designing the algorithm's parameters. Section \ref{sec:DUE} reports numerical experiments in solving dynamic user equilibria in standard test instances, and compares out method with the projection-based algorithm described in \cite{FriHanLiuSze13,HanEveFri19,FriHan19}. 
\section{Preliminaries}
\label{sec:prelims}
We follow the standard notation as in \cite{BauCom16}. Let $\N:=\{1,2,\ldots\}$ be the set of positive integers and $\N_{0}:=\{0\}\cup\N$ the set of nonnegative integers. Let $\scrH$ be a real separable Hilbert space with inner product $\inner{x,y}$ and induced norm $\norm{x}:=\sqrt{\inner{x,x}}$. A sequence $(x_{n})_{n\in\N}$ converges strongly to a point $x\in\scrH$ if $\lim_{n\to\infty}\norm{x_{n}-x}=0$. A sequence $(x_{n})_{n\in\N}$ converges weakly to a point $x\in\scrH$ if, for every $u\in\scrH$, $\inner{x_{n},u}\to \inner{x,u}$; in symbols $x_{n}\wlim x$. 

Let $\scrX\subseteq\scrH$ be a closed convex nonempty subset. Define the normal cone mapping by $
\NC_{\scrX}(x):=\{u\in\scrH\vert \inner{u,y-x}\leq 0\quad\forall y\in\scrX\}$ if $x\in\scrX$, and $\NC_{\scrX}(x)=\emptyset$ otherwise. The Euclidean projector onto $\scrX$ is defined as $P_{\scrX}(x):= \argmin_{y\in\scrX}\frac{1}{2}\norm{y-x}^{2}$. 
It is well known that $P_\scrX$ is nonexpansive and the following property, taken from \cite{GR1984}, hold.

\begin{lemma}\label{b7}
Let $\scrX$ be a nonempty closed convex subset of a real Hilbert space $\scrH.$ Given $x\in \scrH$ and $z\in \scrX$. Then
\begin{equation}
z=P_\scrX (x) \Longleftrightarrow \inner{ x-z,z-y} \geq 0 \qquad \forall y\in \scrX.
\end{equation}
\end{lemma}


\begin{definition}
A mapping $F:\scrH\to\scrH$ is \emph{pseudo-monotone} on $\scrX$ if for all $x,y\in\scrX$ it holds 
\begin{equation}\label{eq:PM}
\inner{F(x),y-x}\geq 0\Rightarrow \inner{F(y),y-x}\geq 0
\end{equation}
The mapping $F:\scrH\to\scrH$ is \emph{monotone} on $\scrX$ if for all $x,y\in\scrX$ it holds 
\begin{equation}\label{eq:M}
\inner{F(x)-F(y),x-y}\geq 0.
\end{equation}
\end{definition}

Clearly, pseudo-monotonicity is a weakened monotonicity assumption providing enough structure to derive provably strongly convergent algorithms. In particular, if $F=\nabla f$ is the gradient of a differentiable real-valued function $f:\scrH\to\R$, pseudo-monotonicty coincides with \emph{pseudo-convexity} of the function $f$. Pseudo-convexity is the classical assumption involved in existence proofs of economic equilibria and Nash equilibria in games with continuous action spaces \cite{FacPan03}.

In the sequel, we use the following simple relations:\\
For each $x, y, z \in \scrH$ and for all $\alpha, \beta, \gamma \in [0, 1]$ with $\alpha + \beta + \gamma = 1$, we have
\begin{align}\label{11x}
\norm{x+y}^2 & \leq \norm{x}^2+2\inner{y,x+y},\text{ and }\\
\label{eee}
\norm{\alpha x+\beta y+\gamma z}^2 &= \alpha\norm{x}^2 + \beta\norm{y}^2 + \gamma\norm{z}^2- \alpha \beta \norm{x-y}^2 \\
& \quad - \alpha \gamma \norm{x-z}^2  -  \beta \gamma\norm{y-z}^2.\nonumber
\end{align}
The next technical lemma, due to Xu \cite{Xu02}, will be key in our convergence analysis.

\begin{lemma}\label{Xu} 
Let $(a_n)_{n\in\N_{0}}$ be sequence of nonnegative real numbers such that:
\begin{equation}
a_{n+1}\leq (1-\alpha_n)a_n+\alpha_n b_n,
\end{equation}
where $(\alpha_n)_{n\in \N_{0}}\subset (0,1)$ and $(b_n)_{n\in\N_{0}}$ is a sequence such that
\begin{itemize}	
\item[(a)] $\sum_{n=0}^\infty \alpha_n=\infty$, and 
\item[(b)] $\limsup_{n\to\infty}b_n\le 0.$
\end{itemize}
Then $\lim_{n\to\infty}a_n=0.$
\end{lemma}
\section{A Strongly Convergent Algorithm for Pseudo-monotone VIs}
\label{sec:VI}
We are given a mapping $F:\scrH\to\scrH$, satisfying the following assumptions:
\begin{assumption}
\label{ass:Lipschitz}
$F:\scrH\to\scrH$ is Lipschitz continuous with Lipschitz constant $L>0$, and sequentially weak-to-weak continuous on bounded subsets of $\scrH$.
\end{assumption}
Recall that weak-to-weak continuity requires that for every weakly converging sequence $x_{n}\rightharpoonup x$, it follows that $F(x_{n})\wlim F(x)$ \citep{BauCom16}. In terms of regularity, we also rely on the following mild monotonicity assumption on the map $F$:
\begin{assumption}
\label{ass:PM}
$F:\scrH\to\scrH$ is \emph{pseudomonotone} on $\scrX$: For all $x,y\in\scrX$ it holds 
\begin{equation}\label{eq:PM}
\inner{F(x),y-x}\geq 0\Rightarrow \inner{F(y),y-x}\geq 0
\end{equation}
\end{assumption}
Our objective is to solve the Hilbert-space valued variational inequality $\VI(\scrX,F)$:
\begin{equation}\label{eq:VI}
\text{ find } x^\ast \in\scrH \text{ such that } \inner{F(x^{\ast}),x-x^{\ast}}\geq 0\qquad \forall x\in\scrX. 
\end{equation} 
\begin{assumption}\label{ass:Solution}
Let $\scrX_{\ast}\subset\scrX$ denote the set of solutions to $\VI(\scrX,F)$. Then $\scrX_{\ast}$ is a nonempty, closed and convex set. 
\end{assumption}
For standard conditions guaranteeing existence of solutions to $\VI(\scrX,F)$ we refer the reader to \cite{AusTebBook06}.

\subsection{Algorithmic Setting}
In this section we present two strongly convergent numerical schemes for solving $\VI(\scrX,F)$ under Assumptions \ref{ass:Lipschitz}-\ref{ass:Solution}. The building block of our construction is the classical forward-backward-forward algorithm proposed by Tseng \cite{Tse00}, in the context of solving monotone inclusions. As is well known, the advantage of Tseng's splitting technique is that it allows us to treat monotone inclusions for finding zeroes of the operator $A+B$, where $A:\scrH\to2^{\scrH}$ and $B:\scrH\to\scrH$ are both maximally monotone and $B$ is $L$-Lipschitz. Compared to the celebrated forward-backward splitting, Tseng's method does not require cocoercivity of the single-valued operator $B$. When applied to variational inequalities, the main advantage of the forward-backward-forward method is that it requires only a single projection step at each iteration, which makes the algorithm much more efficient in practice relative to its close competitor the extragradient method of \cite{Kor76}.\footnote{See \cite{BotCseVuo18,StaMerIFAC19,BotMerStaVuo19} for an in-depth discussion in stochastic and deterministic variational inequality problems.} We first study a non-adaptive version of our strongly convergent forward-backward-forward algorithm (Algorithm \ref{alg:FBF-VI}). This scheme iteratively constructs a sequence $(x^{k},r^{k},z^{k})_{k\in\N_{0}}\subset\scrH\times\scrH\times\scrX$, where $z^{k}$ and $r^{k}$ are just the classical forward-backward-forward iterations. If we would run the scheme only with these two iterative steps, the best we can hope for is weak convergence of the iterates under the common hypothesis that the map $F$ is monotone. The innovative element of our scheme is the additional extrapolation step generating $x^{k+1}$, which will be enforcing strong convergence of the trajectories to a minimum norm solution of $\VI(\scrX,F)$. We would like to point out that this modification of the forward-backward-forward scheme is much simpler than the one presented in \cite{ThoChoVin19} since no hyperplane projection subroutine is involved in our construction. In view of our objective to develop numerical methods with cheap iterations, this is notable feature of our algorithmic approach. The main theoretical result of this paper reads then as follows:

\begin{algorithm}[t]
\caption{FBF for $\VI(F,\scrX)$.}\label{alg:FBF-VI}
\SetAlgoLined
\KwData{step-size sequence $\gamma\in(0,1/L)$, parameters $(\alpha_{k})_{k\in\N_{0}},(\beta_{k})_{k\in\N_{0}}\subset(0,1)$, Map $F:\scrH\to\scrH$.
}
\KwResult{Minimal norm solution $x^{\ast}\in\scrX_{\ast}$ of $\VI(F,\scrX)$.}
Initial point $x^{0}\in\scrX$\;
 \While{$k=0,1,\ldots,k_{\max}$}{
  obtain $x^{k}$\; 
    \eIf{Stopping condition not satisfied}{
Compute $z^{k}=P_{\scrX}[x^{k}-\gamma F(x^{k})]$\; 
Compute $r^{k}=z^{k}+\gamma(F(x^{k})-F(z^{k}))$\;
Update $x^{k+1}=(1-\alpha_{k}-\beta_{k})x^{k}+\beta_{k}r^{k}.$
}{
   Stop and report $x^{k}$ as the solution\;
  }
 }
\end{algorithm}

\begin{theorem}
\label{th:SC}
Let $(\alpha_k)_{k\in\N_{0}}$ and $(\beta_k)_{k\in\N_{0}}$ be two real sequences in $(0,1)$, such that $(\beta_k)_{k\in\N_{0}} \subset (\alpha,1-\alpha_k)$ for some $\alpha>0$, and
\begin{equation}\label{Cond}
\lim_{k\to\infty}\alpha_k=0, \sum_{k=1}^{\infty}\alpha_k=\infty.
\end{equation}
Then the sequence $(x^{k})_{k\in\N_{0}}$ generated by Algorithm \ref{alg:FBF-VI} converges strongly to  $p\in \scrX_{\ast}$, where $p=\argmin\{\|z\|: z\in \scrX_{\ast}\}$.
\end{theorem}

Beside excellent convergence properties and computationally cheap iterations, Algorithm \ref{alg:FBF-VI} requires knowledge of the Lipschitz constant of the map $F$. In practice we usually have no information about such a global quantity, making the applicability of Algorithm \ref{alg:FBF-VI} questionable. Fortunately, we can circumvent this annoying strong assumption by constructing a simple adaptive step-size policy relying on evaluations of the function $F$ only, without requesting explicit knowledge of the Lipschitz constant. Specifically, let us consider a sequence $(\gamma_{k})_{k\in\N_{0}}$, defined recursively by
\begin{equation}\label{eq:adaptiveStep}
\gamma_{k+1}:= \left\{\begin{array}{ll} 
\min\left\{\frac{\rho\norm{z^{k}-x^k}}{\norm{F(z^k)-F(x^k)}}, \gamma_k\right\} &  \text{ if }F(z^k)-F(x^k) \neq 0,\\
\gamma_k  & \text{otherwise.}
\end{array}\right.
\end{equation} 
The parameters $\rho \in (0,1)$ and $\gamma_0$ are chosen at the beginning of the scheme by the user. It is clear that $(\gamma_k)_{k\in\N_{0}}$ is non-increasing and bounded from above by $\min\left\{\gamma_0, \frac{\rho}{L}\right\}$. This implies that the sequence $(\gamma_{k})_{k\in\N_{0}}$ has a limit point not smaller than $\left\{\gamma_0, \frac{\rho}{L}\right\}$. Replacing in Algorithm \ref{alg:FBF-VI} the constant step-size $\gamma$ by the sequence $(\gamma_{k})_{k\in\N_{0}}$, leads us directly an adaptive forward-backward-forward scheme, precisely defined in Algorithm \ref{alg:FBF-adapt}. 

\begin{algorithm}[t]
\caption{FBF for $\VI(F,\scrX)$ adaptive step-size}\label{alg:FBF-adapt}
\SetAlgoLined
\KwData{Initial step-size $\gamma_0 > 0$, parameters $\rho \in (0,1)$,  $(\alpha_{k})_{k\in\N_{0}},(\beta_{k})_{k\in\N_{0}}\subset(0,1)$\;
Map $F:\scrH\to\scrH$.
}
\KwResult{Minimal norm solution $x^{\ast}\in\scrX_{\ast}$ of $\VI(\scrX,F)$.}
Initial point $x^{0}\in\scrX$\;
 \While{$k=0,1,\ldots,k_{\max}$}{
  obtain $x^{k}$\; 
    \eIf{Stopping condition not satisfied}{
    Compute $z^{k}=P_{\scrX}[x^{k}-\gamma_k F(x^{k})]$\;
    Compute $r^{k}=z^{k}+\gamma_k(F(x^{k})-F(z^{k}))$\; 
Update $x^{k+1}=(1-\alpha_{k}-\beta_{k})x^{k}+\beta_{k}r^{k}$\; 
Update new step-size $\gamma_{k+1}$ by \eqref{eq:adaptiveStep}.
}{
   Stop and report $x^{k}$ as the solution.
  }
 }
\end{algorithm}

\begin{theorem}
\label{th:SCadaptive}
Let $(\alpha_k)_{k\in\N_{0}}$ and $(\beta_{k})_{k\in\N_{0}}$ be two real sequences in $(0,1)$, satisfying the same conditions as in Theorem \ref{th:SC}. Let $(\gamma_{k})_{k\in\N_{0}}$ be designed by the adaptive rule \eqref{eq:adaptiveStep}. Then the sequence $(x^{k})_{k\in\N_{0}}$ generated by Algorithm \ref{alg:FBF-VI} converges strongly to $p=\argmin\{\|z\|: z\in \scrX_{\ast}\}$.
\end{theorem}
The proof of this Theorem only requires a simple twist of the proof of Theorem \ref{th:SC}, and is given at the end of the next Section.

\subsection{Convergence Analysis}
This section is devoted to the proof of Theorem \ref{th:SC} and Theorem \ref{th:SCadaptive}. The proof of these two main results require a series of technical auxiliary results which are collected here. The reader interested in the application to DUE can skip this section and go directly to Section \ref{sec:DUE}.

As a first step in our convergence analysis, we need the following, admittedly quite classical, fundamental recursion.
\begin{lemma}\label{classic1}
Let $x^\ast \in \scrX_{\ast}$ be an arbitrary solution of $\VI(\scrX,F)$. Then, for all $k\geq 0$, we have 
\begin{equation}\label{mainIne}
\norm{r^{k} - x^\ast}^2 \leq \norm{x^{k} - x^\ast}^2 - \left(1 - (\gamma L)^2 \right) \norm{x^{k} - z^{k}}^2.
\end{equation}
\end{lemma} 
{\it Proof.} 
Pick $x^\ast \in \scrX_{\ast}$ arbitrary, and $k \geq 0$ be a fixed iteration counter. Since $z^{k} \in\scrX$, we have by definition
$$
\langle F(x^\ast), z^{k} - x^{\ast} \rangle \geq 0.
$$
From pseudo-monotonicity of $F$ (Assumption \ref{ass:PM}), it follows that 
\begin{equation}\label{+}
\inner{F(z^{k}), z^{k} - x^{\ast}} \geq 0.
\end{equation}
On the other hand, since $z^{k} = P_{\scrX}[x^{k} - \gamma F(x^{k})]$, Lemma \ref{b7} gives
\begin{equation}\label{++}
\inner{x^\ast - z^{k}, z^{k} - x^{k} + \gamma F(x^{k})} \geq 0.
\end{equation}
Multiplying both sides of \eqref{+} by $\gamma > 0$, and adding the resulting inequality to \eqref{++}, we arrive at the bound
\begin{equation}
\inner{x^\ast - z^{k}, z^{k} - x^{k} + \gamma F(x^{k}) - \gamma F(z^{k})} \geq 0.
\end{equation}
Equivalently, 
\begin{equation}
\inner{x^\ast - z^{k}, r^{k} - x^{k}} \geq 0. 
\end{equation}
From this, we can deduce that 
\begin{align}
\inner{r^{k} - x^\ast, r^{k} - x^{k}} &\leq \inner{r^{k} - z^{k}, r^{k} - x^{k}} \nonumber \\
&= \norm{r^{k} - x^{k}}^2 + \inner{x^{k} - z^{k}, r^{k} - x^{k}} \nonumber \\
&= \norm{r^{k} - x^{k}}^2 + \inner{x^{k} - z^{k}, z^{k} + \gamma (F(x^{k}) - F(z^{k})) - x^{k}} \nonumber \\
&= \norm{r^{k} - x^{k}}^2 - \norm{z^{k} - x^{k}}^2 + \gamma \inner{x^{k} - z^{k}, F(x^{k}) - F(z^{k})}.  \label{star}
\end{align}
Recall the elementary Pythagoras identity 
\begin{equation}\label{starstar}
\norm{r^{k} - x^\ast}^2 - \norm{x^{k} - x^\ast}^2 + \norm{r^{k} - x^{k}}^2 = 2 \inner{r^{k} - x^\ast, r^{k} - x^{k}}. 
\end{equation} 
Combining \eqref{star} and \eqref{starstar}, we obtain 
\begin{align}
\nonumber
\norm{r^{k}-  x^\ast}^2 &\leq \norm{x^{k} - x^\ast}^2 + \norm{r^{k} - x^{k}}^2 - 2 \norm{z^{k} - x^{k}}^2 \\
&+ 2 \gamma \inner{x^{k} - z^{k}, F(x^{k}) - F(z^{k})}. 
\label{pentagramm}
\end{align}
Using that $F$ is $L$-Lipschitz yields
\begin{align}
\norm{r^{k}- x^{k}}^2 
&= \norm{\gamma (F(x^{k}) - F(z^{k})) + z^{k} - x^{k}}^2 \nonumber \\
&= \norm{z^{k} - x^{k}}^2 + 2 \gamma \inner{z^{k} - x^{k}, F(x^{k}) - F(z^{k})} + \gamma^2 \norm{F(x^{k}) - F(z^{k})}^2 \nonumber \\
&\leq \norm{z^{k}- x^{k}}^2 + 2 \gamma \inner{z^{k} - x^{k}, F(x^{k}) - F(z^{k})} + (\gamma L)^2 \norm{x^{k} - z^{k}}^2.
\label{pentagramm1} 
\end{align}
Finally, combining \eqref{pentagramm} with \eqref{pentagramm1}, we obtain the desired inequality  
\begin{align}
\norm{r^{k} - x^\ast}^2 &\leq \norm{x^{k} - x^\ast}^2 - \norm{z^{k} - x^{k}}^2 + (L\gamma)^2 \norm{x^{k} - z^{k}}^2 \\
&= \norm{x^{k} - x^\ast}^2 - \left(1 - (L\gamma)^2\right) \norm{z^{k} - x^{k}}^2
\end{align}
\qed

Next, we establish boundedness of the produced trajectory $(x^{k})_{k\in\N_{0}}$. 

\begin{lemma}
The sequence $(x^k)_{k\in\N_{0}}$ generated by Algorithm \ref{alg:FBF-VI} is bounded.
\end{lemma}
{\it Proof.} 
Thanks to Lemma \ref{classic1} and $\gamma\in(0,1/L)$, we have for every $x^\ast \in \scrX_{\ast}$ 
\begin{equation}\label{qvx2}
\norm{r^k-x^*}\leq\norm{x^k-x^*} \quad \forall k\geq 0.
\end{equation}
By definition of the iterate $x^{k+1}$, the triangle inequality gives us 
\begin{align}\label{qt10}
\norm{x^{k+1}-x^*}&=\norm{(1-\alpha_k-\beta_k)x^k+\beta_k r^k- x^{\ast}}\\
&=\norm{(1-\alpha_k-\beta_k)(x^k-x^*)+\beta_k(r^k-x^{\ast})-\alpha_k x^*}\notag\\
&\leq \norm{(1-\alpha_k-\beta_k)(x^k-x^*)+\beta_k(r^k-x^*)}+\alpha_k \norm{x^{\ast}}.
\end{align}
From \eqref{qvx2}, we obtain for all $k \geq 0$ that
\begin{align}
\norm{(1-&\alpha_k-\beta_k)(x^k-x^*)+\beta_k(r^k-x^*)}^2\\
=&(1-\alpha_k-\beta_k)^2\norm{x^k-x^*}^2+2(1-\alpha_k-\beta_k)\beta_k \inner{x^k-x^*,r^k-x^* }+\beta^2_k\norm{r^k-x^{\ast}}^2\\
\leq& (1-\alpha_k-\beta_k)^2 \norm{x^k-x^*}^2+2(1-\alpha_k-\beta_k)\beta_k  \norm{r^k-x^*}\cdot \norm{x^k-x^*} +\beta^2_k \norm{r^k-x^*}^2\\
\leq & (1-\alpha_k-\beta_k)^2 \norm{x^k-x^*}^2+2(1-\alpha_k-\beta_k)\beta_k \norm{x^k-x^*}^2 +\beta^2_k\norm{x^k-x^*}^2\\
=& (1-\alpha_k)^2 \norm{x^k-x^*}^2.
\end{align}
This implies 
\begin{equation}\label{qt9}
\norm{(1-\alpha_k-\beta_k)(x^k-x^*)+\beta_k(r^k-x^{\ast})} \leq (1-\alpha_k) \norm{x^k-x^{\ast}} \quad \forall k \geq 0.
\end{equation}
Combining (\ref{qt10}) and (\ref{qt9}), we get by induction
\begin{align}
\norm{x^{k+1}-x^*}&\leq (1-\alpha_k) \norm{x^k-x^*}+\alpha_k\norm{x^*} \\
&\leq \max\{\norm{x^k-x^*},\norm{x^*} \}\\
&\vdots\\
&\leq\max\{\norm{x_{0}-x^*},\norm{x^*}\}.
\end{align}
Hence, we conclude that the sequence $(x^k)_{k\in\N_{0}}$ is bounded, and so is $(r^k)_{k\in\N_{0}}$.
\qed

The next lemma is key to the proof of the main result of this paper. While at first sight it looks very similar to typical bounds obtained in the setting of quasi-Fej\'{e}r iterations, it provides us the necessary structure to deduce strong convergence of the iterates via Lemma \ref{Xu}. 
 
\begin{lemma} \label{Lem9}
Let $(\alpha_{k})_{k\in\N_{0}}$ and $(\beta_{k})_{k\in\N_{0}}$ be two sequences satisfying $\beta_k < 1 - \alpha_k$ for all $k\in\N_{0}$. Then, for all  $x^\ast \in \scrX_{\ast}$ and $k\geq 0$, it holds
\begin{equation}
\norm{x^{k+1}-x^*}^2\leq (1-\alpha_k)\norm{x^k-x^*}^2+\alpha_k\left [2\beta_k\norm{x^k-r^k}\cdot\norm{x^{k+1}-x^{\ast}}+2\inner{x^*,x^*-x^{k+1}} \right].
\end{equation}
\end{lemma}
{\it Proof.} 
Using (\ref{eee}), we have
\begin{align}
\norm{x^{k+1}-x^*}^2=&\norm{(1-\alpha_k-\beta_k)x^k+\beta_k r^k- x^*}^2\notag\\
=&\norm{(1-\alpha_k-\beta_k)(x^k-x^*)+\beta_k (r^k- x^*)+\alpha_k (-x^*)}^2\notag\\
=&(1-\alpha_k-\beta_k)\norm{x^k-x^*}^2+\beta_k \norm{r^k- x^*}^2+\alpha_k \norm{x^*}^2 \notag\\
&-\beta_k (1-\alpha_k-\beta_k)\norm{x^k-r^k}^2-\alpha_k(1-\alpha_k-\beta_k) \norm{x^k}^2-\alpha_k \beta_k \norm{r^k}^2\notag\\
\leq& (1-\alpha_k-\beta_k)\norm{x^k-x^*}^2+\beta_k \norm{r^k- x^*}^2+\alpha_k\norm{x^*}^2.
\label{qxv5}
\end{align}
Together with Lemma \ref{classic1}, this implies
\begin{align}
\norm{x^{k+1}-x^*}^2& \leq (1-\alpha_k-\beta_k)\norm{x^k-x^*}^2+\beta_k \norm{x^k- x^*}^2\notag \\
&\quad -\beta_k \left(1-(\gamma L)^2  \right) \norm{x^k-z^k}^2+\alpha_k \norm{x^*}^2\notag\\
&=(1-\alpha_k)\norm{x^k-x^*}^2-\beta_k \left(1-(\gamma L)^2  \right) \norm{x^k-z^k}^2+\alpha_k \norm{x^*}^2\notag\\
&\leq \norm{x^k-x^*}^2-\beta_k \left(1-(\gamma L)^2  \right) \norm{x^k-z^k}^2+\alpha_k \norm{x^*}^2.
\end{align}
Therefore, 
\begin{equation} \label{V1}
\beta_k \left(1-(\gamma L)^2  \right) \norm{x^k-z^k}^2\leq \norm{x^k-x^*}^2-\norm{x^{k+1}-x^*}^2 +\alpha_k \norm{x^*}^2.
\end{equation}
Setting $t^k=(1-\beta_k)x^k+\beta_k r^k$ we obtain
\begin{align}\label{qxv9}
\norm{t^k-x^*}=&\norm{(1-\beta_k)(x^k-x^*)+\beta_k(r^k-x^*)}\\
&\leq (1-\beta_k)\norm{x^k-x^*}+\beta_k\norm{r^k-x^*}\notag\\
&\leq (1-\beta_k)\norm{x^k-x^*}+\beta_k\norm{x^k-x^*}  \\
=&\norm{x^k-x^*},
\end{align}
and 
\begin{equation}\label{qxv10}
\norm{t^k-x^k}=\beta_k \norm{x^k-r^k}.
\end{equation}
Combining (\ref{qxv9}) with (\ref{qxv10}), using \eqref{11x}, we get
\begin{align}
\norm{x^{k+1}-x^*}^2&=\norm{(1-\alpha_k-\beta_k)x^k+\beta_k r^k-x^*}^2\\
&=\norm{(1-\beta_k)x^k+\beta_k r^k-\alpha_k x^k-x^*}^2\\
&=\norm{(1-\alpha_k)(t^k-x^*)-\alpha_k(x^k-t^k)-\alpha_k x^*}^2\\
& \leq (1-\alpha_k)^2\norm{t^k-x^*}^2 -2\inner{\alpha_k(x^k-t^k)+\alpha_k x^*,x^{k+1}-x^*} \\
& = (1-\alpha_k)^2\norm{t^k-x^*}^2 +2 \alpha_k\inner{x^k-t^k,x^*-x^{k+1}}+2\alpha_k \inner{x^*,x^*-x^{k+1}}\\
& \leq (1-\alpha_k)\norm{t^k-x^*}^2 + 2\alpha_k \norm{x^k-t^k}\cdot \norm{x^{k+1}-x^*}+2\alpha_k \inner{x^*,x^*-x^{k+1}}\\
& \leq (1-\alpha_k)\norm{x^k-x^*}^2 + \alpha_k \left [2\beta_k\norm{x^k-r^k}\cdot\norm{x^{k+1}-x^*}+2\inner{x^*,x^*-x^{k+1}}\right].
\end{align}
Observe that the assumption $\alpha_{k}\in(0,1)$ has been used here as well.
\qed

The following fundamental result relies heavily on the pseudo-monotonicity and weak continuity of $F$. 
\begin{lemma} \label{Lem10}
Assume there exists a subsequence $(x^{k_j})_{j\in\N}$ of $(x^k)_{k\in\N_{0}}$ such that $(x^{k_j})_{j\in\N}$ converges weakly to $\hat{x}$. Let $(z^{k_{j}})_{j\in\N}$ the corresponding subsequence of $(z^{k})_{k\in\N_{0}}$. If $\lim_{j \to \infty} \norm{x^{k_j}-z^{k_j}}=0$, then $\hat{x} \in \scrX_{\ast}$. 
\end{lemma}
{\it Proof.} 

Let $(x^{k_{j}})_{j\in\N}$ be a converging subsequence with weak limit $\hat{x}$. Since $\lim_{j \to \infty}\norm{x^{k_j}-z^{k_j}}=0$, $(z^{k_j})_{j\in\N}$ also converges weakly to $\hat{x}$. By definition, $(z^{k_j})_{j\in\N} \subset \scrX$ and $\scrX$ is weakly closed. Hence, $\hat{x} \in \scrX$, and we have to prove that $\hat{x} \in \scrX_{\ast}$. Indeed, since for all $j\in\N$, 
\[
z^{k_j} =P_{\scrX}(x^{k_j}-\gamma F(x^{k_j})),
\]
we have
\[
\inner{x^{k_j}-\gamma F(x^{k_j})-z^{k_j}, y-z^{k_j}} \leq 0, \qquad \forall y \in \scrX, 
\]
or equivalently,
\[
\frac{1}{\gamma}\inner{x^{k_j}-z^{k_j}, y-z^{k_j}}\leq\inner{F(x^{k_j}),y-z^{k_j}} \qquad \forall y \in \scrX. 
\]
This implies that
\begin{equation}\label{ine-1}
\frac{1}{\gamma}\inner{x^{k_j}-z^{k_j}, y-z^{k_j}}\leq\inner{F(x^{k_j})-F(z^{k_j}),y-z^{k_j}}+\inner{F(z^{k_j}),y-z^{k_j}}\qquad \forall y \in \scrX. 
\end{equation}
Fixing $y \in \scrX$ and letting $j \to +\infty$ in the last inequality, remembering that  $\lim_{j \to \infty}\norm{x^{k_j}-z^{k_j}}=0$ and $\lim_{j \to \infty} \norm{F(x^{k_j})-F(z^{k_j})}=0$ (by weak-to-weak continuity of $F$), we have
\begin{equation}\label{ine1}
\liminf_{j\to \infty}\inner{F(z^{k_j}),y-z^{k_j}}\geq 0\qquad\forall y\in\scrX.
\end{equation}
Next, choose a sequence $(\epsilon_{j})_{j\in\N}\subset(0,\infty)$ with $\epsilon_{j}\downarrow 0$.  Construct a sequence $(N_{j})_{j\in\N}\subset\N$ such that 
	\begin{equation}\label{ine2}
	\inner{F(z^{k_{i}}),y-z^{k_{i}} }+ \epsilon_{j}  \geq 0 \quad \forall i\geq N_{j}.
	\end{equation}
The existence of such a sequence follows from \eqref{ine1}. Since  $(\epsilon_{j})_{j\in\N}$ is decreasing, it is easy to see that the sequence $(N_{j})_{j\in\N}$ is increasing. Furthermore, for each $j\geq 1$, $F(z^{N_j})\neq 0$. Setting
\[
u^{N_j}:=\frac{F(z^{N_j})}{\norm{F(z^{N_j})}^2},
\]
we have $\inner{F(z^{N_j}),u^{N_j}}=1$ for each $j\geq 1$. Now we can deduce from \eqref{ine2} that for each $j\in\N$
\[
\inner{F(z^{N_j}),y+\epsilon_{j} u^{N_j} -z^{N_j}}\geq 0.
\]
Since $F$ is pseudo-monotone, this implies that
\begin{equation}\label{ine3}
\inner{F(y+\epsilon_{j} u^{N_j}), y+\epsilon_{j} u^{N_j} -z^{N_j}}\geq 0.
\end{equation}
On the other hand, we have that $(z^{N_j})_{j\in\N}$ converges weakly to $\hat{x}$ when $j \to \infty$. Since $F$ is sequentially weak-to-weak continuous on $\scrH$, wave $F(z^{N_j})\wlim F(\hat{x})$. If $F(\hat{x})=0$, then $\hat{x}\in\scrX_{\ast}$. Hence, let us assume that $F(\hat{x})\neq 0$. Since the norm mapping is sequentially weakly lower semicontinuous, we have 
\[
\norm{F(\hat{x})} \leq \liminf_{j \to \infty}\norm{F(z^{N_j})}.
\]
Since $(z^{N_j})_{j\in\N}\subset(z^{k_j})_{j\in\N}$ and $\epsilon_{j} \downarrow 0$ as $j\to\infty$, we obtain
\[
0 \leq \lim_{j \to \infty} \norm{\epsilon_{j} u^{N_j}} = \lim_{j \to \infty} \frac{\epsilon_{j}}{\norm{F(z^{N_j})}}=0.
\]
Hence, taking the limit as $j\to \infty$ in \eqref{ine3}, we obtain 
\[
\inner{F(y),y-\hat{x}}\geq 0.
\]
This and the pseudomonotonicity and continuity of $F$ imply that $\hat{x} \in \scrX_{\ast} $.
\qed

We are now in the position to prove the main result of this section.\\
{\it Proof of Theorem \ref{th:SC}}
Since $\scrX_{\ast}$ is closed and convex, there exists a unique element $p\in \scrX_{\ast}$ such that $p = P_{\scrX_{\ast}} (0)$. We will show that the sequence $\left(\norm{x^k-p}^2\right)_{k\in\N_{0}}$ converges to zero by considering two possible cases on its long-run behavior. 

{\bf Case 1:} There exists an $k_0\in\N$ such that $\norm{x^{k+1}-p}^2\leq \norm{x^k-p}^2$ for all $k\geq k_0.$ This implies that $\lim_{n\to\infty}\norm{x^k-p}^2$ exists. It follows from \eqref{V1} and \eqref{Cond} that
\begin{equation}\label{v1}
\lim_{k\to\infty} \norm{x^k-z^k}=0.
\end{equation}
We also have
\begin{align}\label{vvn12}
\norm{r^k-x^k}&=\norm{z^k-\gamma (F(z^k)-F(x^k))-x^k} \\
&\leq (1+\gamma L)\norm{x^k-z^k}.
\end{align}
Combining \eqref{v1} and \eqref{vvn12}, we get
\begin{equation}
\lim_{n\to\infty}\norm{r^k-x^k}=0.
\end{equation}
Therefore,
\begin{equation}
\norm{x^{k+1}-x^k}\leq\alpha_k \norm{x^k}+\beta_k\norm{x^k-r^k}\to 0 \text { as } k\to \infty.
\end{equation}
Since $(x^k)_{k\in\N_{0}}$ is bounded, we can without loss of generality assume that there exists a subsequence $(x^{k_j})_{j\in\N}$ such that $x^{k_j}\wlim q$, and
\begin{equation}
\limsup_{k\to\infty}\inner{p,p-x^{k}}=\lim_{j\to\infty}\inner{p,p-x^{k_j}}=\inner{p,p-q}.
\end{equation}
By Lemma \ref{Lem10}, we conclude that $q \in \scrX_{\ast}$. 
Since $p=P_{\scrX_{\ast}}(0)$, we obtain
\begin{equation}
\limsup_{k\to\infty}\inner{ p,p-x^{k}}=\inner{p,p-q}\leq 0.
\end{equation}
Since $\norm{x^{k+1}-x^k}\to 0$, it also must be true that 
 \begin{equation}
 \limsup_{k\to\infty}\inner{p,p-x^{k+1}}\leq 0.
 \end{equation}
 From Lemma \ref{Lem9} and Lemma \ref{Xu}, we finally conclude $\lim_{k\to\infty}\norm{x^{k}-p}^2=0$. That is $x^k\to p.$
	
{\bf Case 2:} Assume that there is no $k_0 \in \mathbb{N}$ such that $\left(\norm{x^{k}-p}\right)_{k\geq k_0}$ is monotonically
decreasing. We follow the technique in \cite{Mai08}; Set $\Gamma_k:=\norm{x^k-p}^2$ for all $k \geq 1$, and let $\tau: \N \to\N$ be a mapping defined for all $k \geq k_0$ (for some $k_0$ large enough) by
\begin{equation}
\tau (k) :=\max\left\{j \in\N\vert j \leq  k, \Gamma_j \leq \Gamma_{j+1}\right\}.
\end{equation}
Hence, $\tau (k)$ is the largest number $j$ in $\{1, . . . , k\}$ such that $\Gamma_j$ increases at $j = \tau (k)$. Note that $\tau (k)$ is well-defined for all sufficiently large $k$. From \cite{Mai08} we deduce that $\left(\tau(k)\right)_{k\in\N}$ is a non-decreasing sequence such that $\tau (k)\to \infty$ as $k \to \infty$, and for all $k\geq k_0$ it holds that
\begin{align}
&0\leq \Gamma_{\tau(k)}\leq \Gamma_{\tau(k)+1},\\
&0 \leq \Gamma_k \le \Gamma_{\tau(k)+1}\label{mas}.
\end{align}
Since $\beta_k\geq  \alpha$ for all $ k\geq 1$, from \eqref{V1} we have
\begin{align}
\alpha \left(1-\gamma_{\tau(k)}^2 L^2  \right)\|x^{\tau(k)}-y^{\tau(k)}\|^2
\leq &\beta_{\tau(k)} \left(1-\gamma_{\tau(k)}^2 L^2  \right)\|x^{\tau(k)}-y^{\tau(k)}\|^2\\
\leq& \norm{x^{\tau(k)}-p}^2-\norm{x^{\tau(k)+1}-p}^2 +\alpha_{\tau(k)} \norm{p}^2\\
\leq &\alpha_{\tau(k)} \norm{p}^2.
\end{align}
Therefore
\[
\lim_{k\to\infty}\norm{x^{\tau(k)}-y^{\tau(k)}}=0 .
\]
As proved in the first case, we have
\begin{equation}
\norm{x^{\tau(k)+1}-x^{\tau(k)}}\to 0,
\end{equation}
and
\begin{equation}
\limsup_{k\to\infty}\inner{p,p-x^{\tau(k)+1}} \leq 0.
\end{equation}
From Lemma \ref{Lem9}  and $\Gamma_{\tau(k)}\le \Gamma_{\tau(k)+1}$ for all $k\geq k_0$, we have
\begin{align}
\|x^{\tau(k)+1}-p\|^2\le& (1-\alpha_{\tau(k)})\|x^{\tau(k)}-p\|^2\\
&+\alpha_{\tau(k)}\left [2\beta_{\tau(k)}\|x^{\tau(k)}-z^{\tau(k)}\|\|x^{\tau(k)+1}-p\|+2\langle p,p-x^{\tau(k)+1}\rangle \right]\\
\leq& (1-\alpha_{\tau(k)})\|x^{\tau(k)+1}-p\|^2\\
&+\alpha_{\tau(k)}[2\beta_{\tau(k)}\|x^{\tau(k)}-z^{\tau(k)}\|\|x^{\tau(k)+1}-p\|+2\langle p,p-x^{\tau(k)+1}\rangle].
\end{align}
This implies that
\begin{equation}
\|x^{\tau(k)+1}-p\|^2\le 2\beta_{\tau(k)}\|x^{\tau(k)}-z^{\tau(k)}\|\|x^{\tau(k)+1}-p\|+2\langle p,p-x^{\tau(k)+1}\rangle,
\end{equation}
which implies that $\limsup_{k\to\infty}\norm{x^{\tau(k)+1}-p}^2\leq 0$; That is 
$$
\lim_{k\to\infty}\norm{x^{\tau(k)+1}-p}=0.
$$
The conclusion follows from \eqref{mas}.
\qed

{\it Proof of Theorem \ref{th:SCadaptive}}

 For the convergence analysis, instead of \eqref{mainIne}, we have 
	\begin{align}
	\|r^{k} - x^\ast\|^2 \leq \|x^{k} - x^\ast\|^2 - \left(1-\frac{\gamma_k^2 \rho^2}{\gamma_{k+1}^2} \right) \|x^{k} - z^{k}\|^2 \qquad\forall k\geq 1.
	\end{align} 
	The rest of the proofs follows analagously to the constant stepsize case, and thus left to the reader. 
\qed
\section{Application to Computing Dynamic User Equilibria}
\label{sec:DUE}
In this section we apply the strongly-convergent forward-backward-forward algorithm to compute dynamic user equilibria in two standard test examples taken from the literature. Our description follows the recent survey \cite{FriHan19}. The numerical examples have been constructed based on the MATLAB package https://github.com/DrKeHan/DTA, documented in \cite{HanEveFri19}.

\subsection{Problem Formulation}
Let $[t_{0},t_{1}]$ be a fixed planning horizon. We are given a connected directed graph $G=(V,A)$ with finite set of vertices $V$, representing traffic intersections (junctions) and arc set $A$, representing road segments. A path $p$ in the graph $G$ is identified with a non-repeating finite sequence of arcs which connect a sequence of different vertices. Hence, an arbitrary path $p$ is identified with the list of edges incident to it, i.e. 
$p=\{a_{1},a_{2},\ldots,a_{m}\}.$ The integer $m=m(p)$ denotes the length of the path $p$. We denote the set of all paths of interest by $\scrP$, and set $\setH:= \R^{\abs{\scrP}}$. We are interested in paths which connect a set of distinguished vertices acting as the \emph{origin-destination} (o/d) pairs in our graph. We are given $N$ distinct o/d pairs denoted as $w_{1},\ldots,w_{N}$, where each $w_{i}=(o_{i},d_{i})\in V$. Call $\scrW:=\{w_{1},\ldots,w_{N}\}$, and the set of paths connecting the o/d pair $w$ is denoted by $\scrP_{w}\subseteq\scrP.$ For each o/d pair $w\in\scrW$ we are given a \emph{demand} $Q_{w}>0$; This represents the number of drivers who have to travel from the origin to the destination described by $w$. For simplicity we assume that this demand is exogenously given. The list $Q=(Q_{w})_{w\in\scrW}$ is often called the \emph{trip table}. In DUE modeling, the single most crucial ingredient is the path delay operator, which maps a given vector of departure rates (path flows) $h$ to a vector of path travel times. We stipulate that path flows are square integrable functions over the planning horizon, so that $h_{p}\in L^{2}([t_{0},t_{1}];\R_{+})$ and $h=(h_{p};p\in\scrP)\in \scrH:= L^{2}([t_{0},t_{1}];\setH)$. To measure the delay of drivers on paths, we introduce the operator $D:\scrH\to\scrH,h\mapsto D(h)$, with the interpretation that $D_{p}(t,h)$ is the path travel time of a driver departing at time $t$ from the origin of path $p$, and following this path throughout. This operator is the result of a dynamic network loading procedure, which is an integrated subroutine in the dynamic traffic assignment problem. See \cite{HanEveFri19} for further information.

On top of path delays, we consider penalty terms of the form $\rho(t+D_{p}(t,h)-T_{A}),$ penalizing all arrival times different from the target time $T_{A}$ (i.e. the usual time of a trip on the o/d. pair $w$). The function $\rho:[-\infty,\infty)\to[0,\infty]$ should be monotonically increasing with $\rho(x)>0$ for $x>0$ and $\rho(x)=0$ for $x\leq 0$. Define the effective delay operator as 
\begin{equation}
\Psi_{p}(t,h):=D_{p}(t,h)+\rho(t+D_{p}(t,h)-T_{A}).
\end{equation}
We thus obtain an operator $\Psi:\scrH\to\scrH$, mapping each profile of path departure rates $h$ to effective delays $\Psi(h)\in\scrH$.

We follow the perceived DUE literature, and stipulate that Wardrop's first principle holds: Users of the network aim to minimize their own travel time, given the departure rates in the system. Thus, a user equilibrium is envisaged, where the delays (interpreted as costs) of all travelers in the same o/d pair are equal, and no traveler can lower his/her costs by unilaterally switching to a different route. To put this behavioral axiom into a mathematical framework, we first formulate the meaning of "minimal costs" in the present Hilbert space setting. Recall the essential infimum of a measurable function $g:[t_{0},t_{1}]\to\R$ as 
$\essinf\{g(t)\vert t\in[t_{0},t_{1}]\}=\sup\left\{x\in\R\vert \Leb(\{s\in[t_{0},t_{1}]:g(s)<x\})=0\right\},$ where $\Leb(\cdot)$ denoted the Lebesgue measure on the real line. Given a profile $h\in\scrH$, define 
\begin{align}
\nu_{p}(h)&:=\essinf\{\Psi_{p}(t,h)\vert t\in[t_{0},t_{1}]\}\qquad\forall p\in\scrP, \text{ and }\\
\nu_{w}(h)&:= \min_{p\in\scrP_{w}}\nu_{p}(h)\qquad\forall w\in\scrW.
\end{align}
On top of minimal costs, we have to restrict the set of departure rates to functions satisfying a basic flow conservation property. Specifically, insisting that all trips are realized, we naturally define the set of feasible flows as
\begin{equation}
\Lambda:=\{f\in\scrH\vert \sum_{p\in\scrP_{w}}\int_{t_{0}}^{t_{1}}f_{p}(t)\dif t=Q_{w}\quad\forall w\in\scrW\}. 
\end{equation}
\begin{definition}
\label{def:DUE-V1}
A profile of departure rates $h^{\ast}\in\scrH$ is a DUE if 
\begin{itemize}
\item[(a)] $h^{\ast}\in\Lambda$, and 
\item[(b)] $h^{\ast}_{p}(t)>0\Rightarrow \Psi_{p}(t,h^{\ast})=\nu_{w}(h^{\ast}).$
\end{itemize}
\end{definition}
In \cite{Frie93} it is observed that the definition of DUE can be formulated equivalently as a variational inequality $\VI(\Lambda,\Psi)$: A flow $h^{\ast}\in\Lambda$ is a DUE if 
\begin{equation}\label{eq:DUE}
\inner{\Psi(h^{\ast}),h-h^{\ast}}\geq 0\qquad\forall h\in\Lambda
\end{equation}
\subsection{A Strongly Convergence Forward-Backward-Forward Scheme for DUE}
Departing from \eqref{eq:DUE}, our aim is to solve the DUE problem by using our strongly convergent forward-backward-forward scheme \ref{alg:FBF-VI}. Adapting this scheme to the usual notation in DUE, we arrive at Algorithm \ref{alg:FBF-DUE}.
\begin{algorithm}[t]
\SetAlgoLined
\KwData{Graph $G=(V,A)$ with o/d pairs $\scrW\subset V\times V$\;
Trip Table $(Q_{w})_{w\in\scrW}$\;
 step-size $\gamma>0$\;
 parameters $(\alpha_{k})_{k\geq 0},(\beta_{k})_{k\geq 0}\subset\R_{+}$
}
\KwResult{An approximate DUE $h^{\ast}$}
Initial path flow $h^{0}\in\scrH$\;
 \While{$k=1,2,\ldots,k_{\max}$}{
  obtain $h^{k}$\; 
  Compute $\epsilon_{k}=\frac{\norm{h^{k+1}-h^{k}}^{2}}{\norm{h^{k}}^{2}}$\;
  \eIf{If $\epsilon_{k}>10^{-4}$}{
   Compute the effective path delays $\Psi_{p}(t,h^{k})$\;
Compute $z^{k}=P_{\Lambda}[h^{k}-\gamma \Psi(h^{k})]$\;
Compute the effective path delays $\Psi_{p}(t,z^{k})$\;
Compute $r^{k}=z^{k}+\gamma(\Psi(h^{k})-\Psi(z^{k}))$\;
Compute $h^{k+1}=(1-\alpha_{k}-\beta_{k})h^{k}+\beta_{k}r^{k}$
}{
   Stop and report $h^{k}=h^{\ast}$ as the solution.
  }
 }
 \caption{Forward-backward-forward algorithm for computing DUE.}\label{alg:FBF-DUE}
\end{algorithm}
Some remarks on the implementation of this algorithm are in order. First, it should be pointed out that Algorithm \ref{alg:FBF-DUE} requires two evaluations of the delay operator $\Psi$. 
As already said, this operator is the outcome of an inner procedure, solving the dynamic network loading part of the model. Dynamic network loading is a separate computational step in the dynamic traffic assignment problem. A very popular formulation of dynamic network loading is the fluid dynamic approximation of traffic flows, known as the Lighthill-Whitham-Richards (LWR) model. We refer the interested reader to \cite{GarPic06} for modeling approaches of the dynamic network loading procedure. In case of the popular LWR model evaluating the delay operator requires solving a coupled system of hyperbolic partial differential equations for the traffic density. It is clear that this procedure is the most costly step in the implementation of Algorithm \ref{alg:FBF-DUE}.\\
Algorithm \ref{alg:FBF-DUE} is, modulo the obvious change in notation, equivalent to Algorithm \ref{alg:FBF-VI} if the delay operator $\Psi$ is Lipschitz continuous and pseudomonotone.

\subsection{Numerical Experiments}
\begin{figure}
    \centering
       \includegraphics[width=1\textwidth]{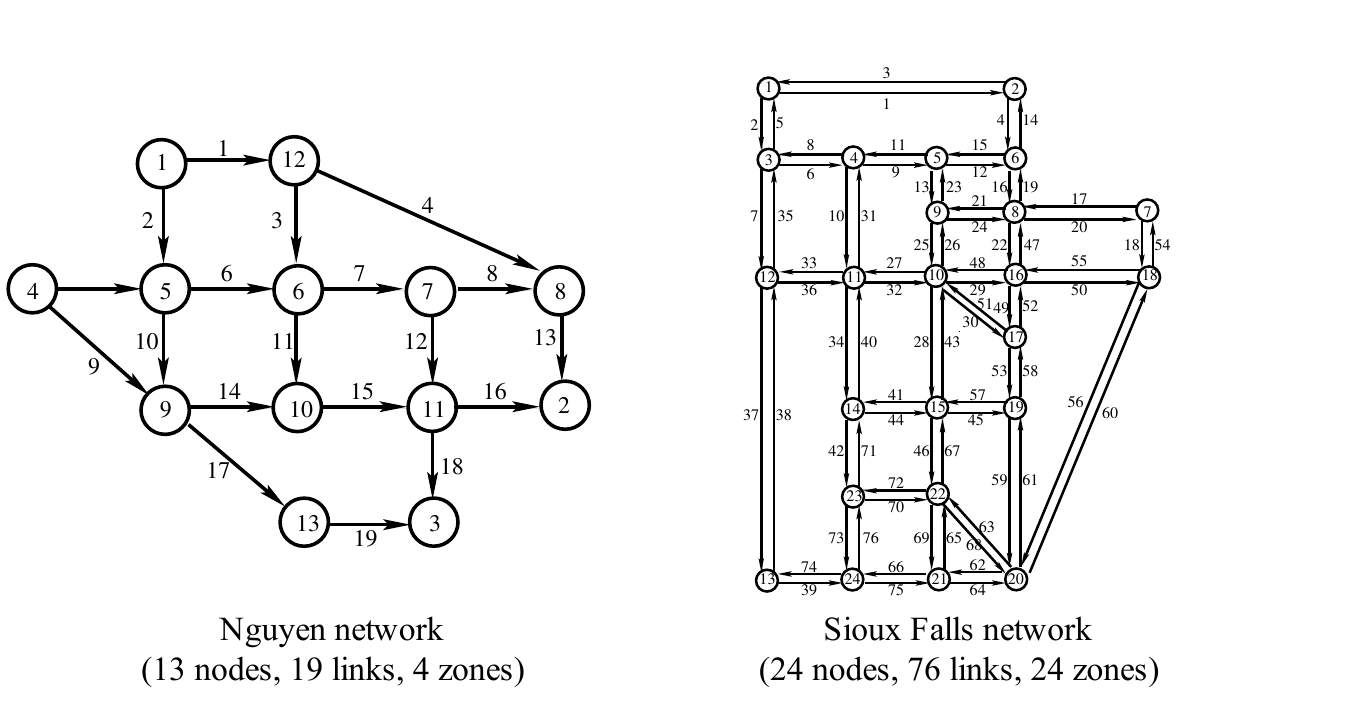}
	\caption{The Nguyen and Sioux Falls network.}
	\label{fig:networks}
	\end{figure}

\begin{figure}
    \centering
        \includegraphics[width=\textwidth]{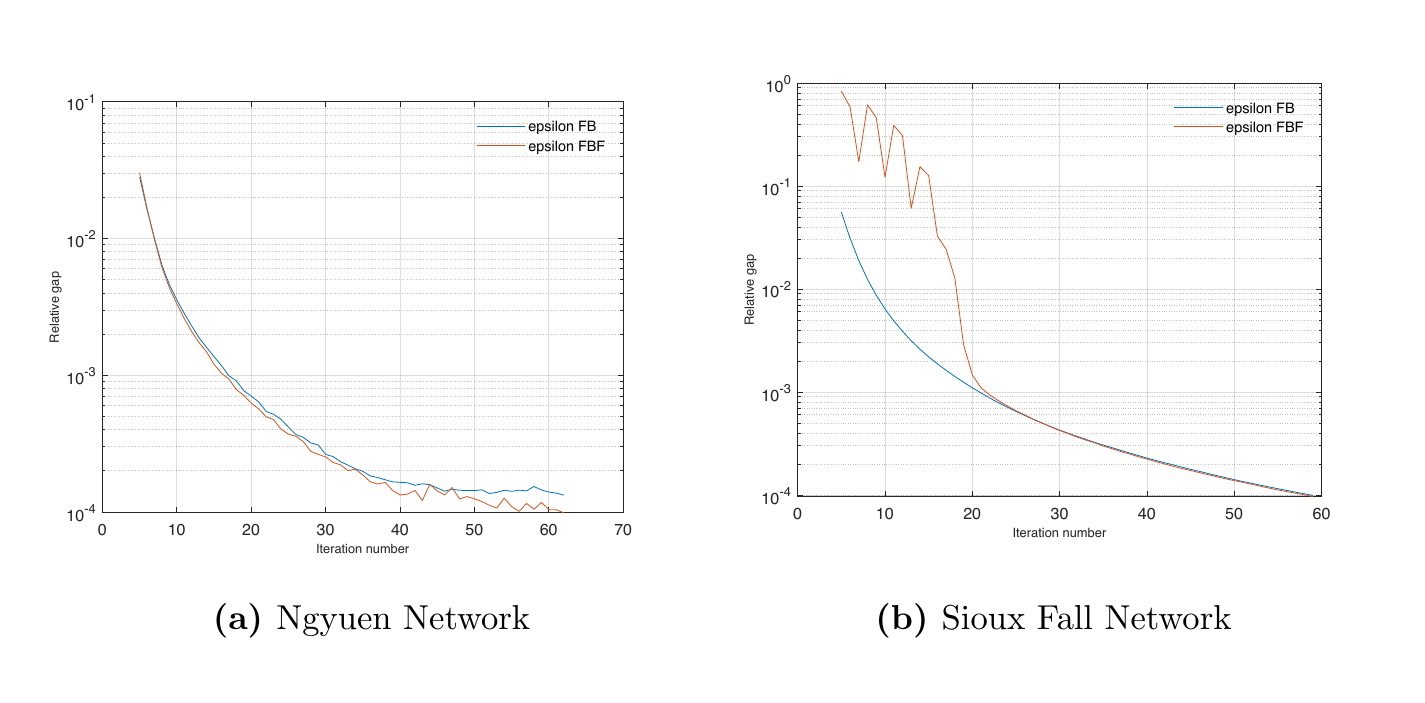}
  \caption{Relative gap \eqref{eq:gap} (called epsilon in the figure) computed under the forward-backward iteration of \cite{HanEveFri19} and Algorithm \ref{alg:FBF-DUE}, using the same parameter values}
   \label{fig:gaps}
\end{figure}

\begin{figure}
    \centering
        \includegraphics[width=\textwidth]{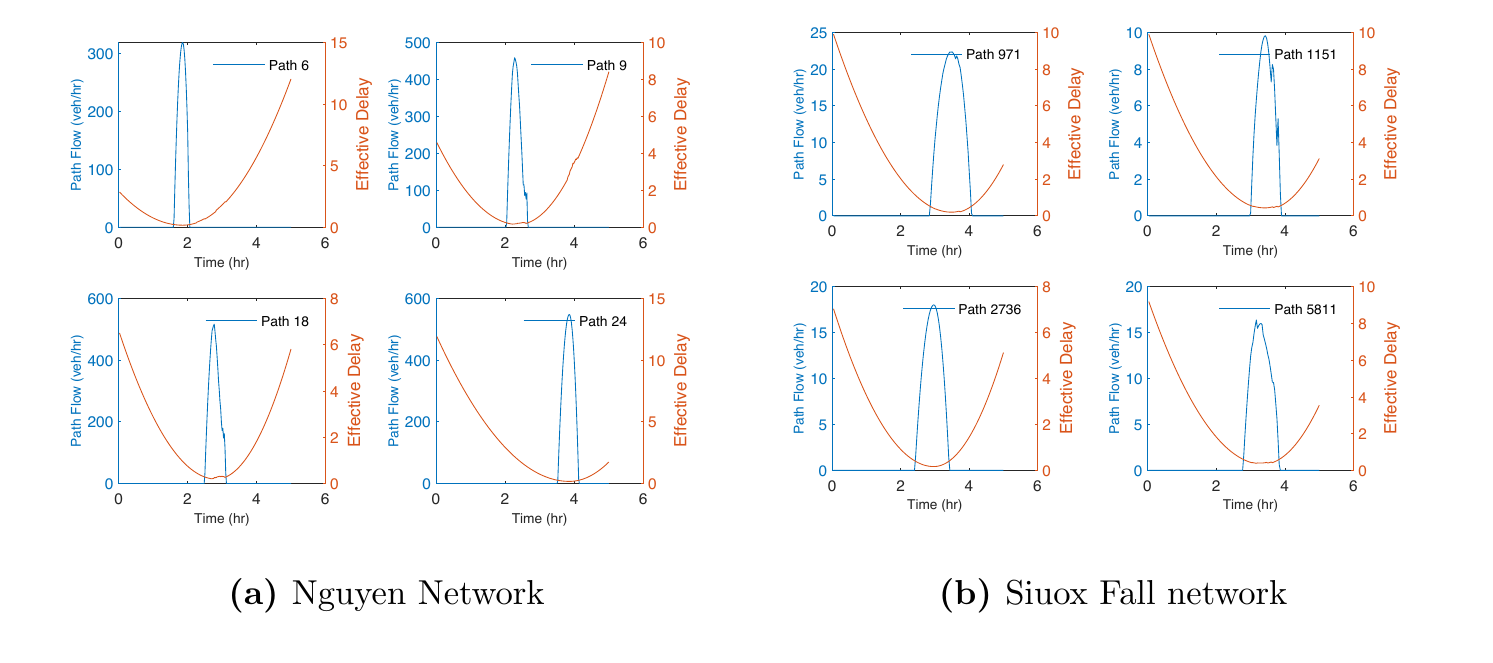}
    \caption{Path departure rates and corresponding effective path delays of selected paths in the DUE solutions.}
    \label{fig:paths}
\end{figure}

\begin{figure}
    \centering
           \includegraphics[width=0.8\textwidth]{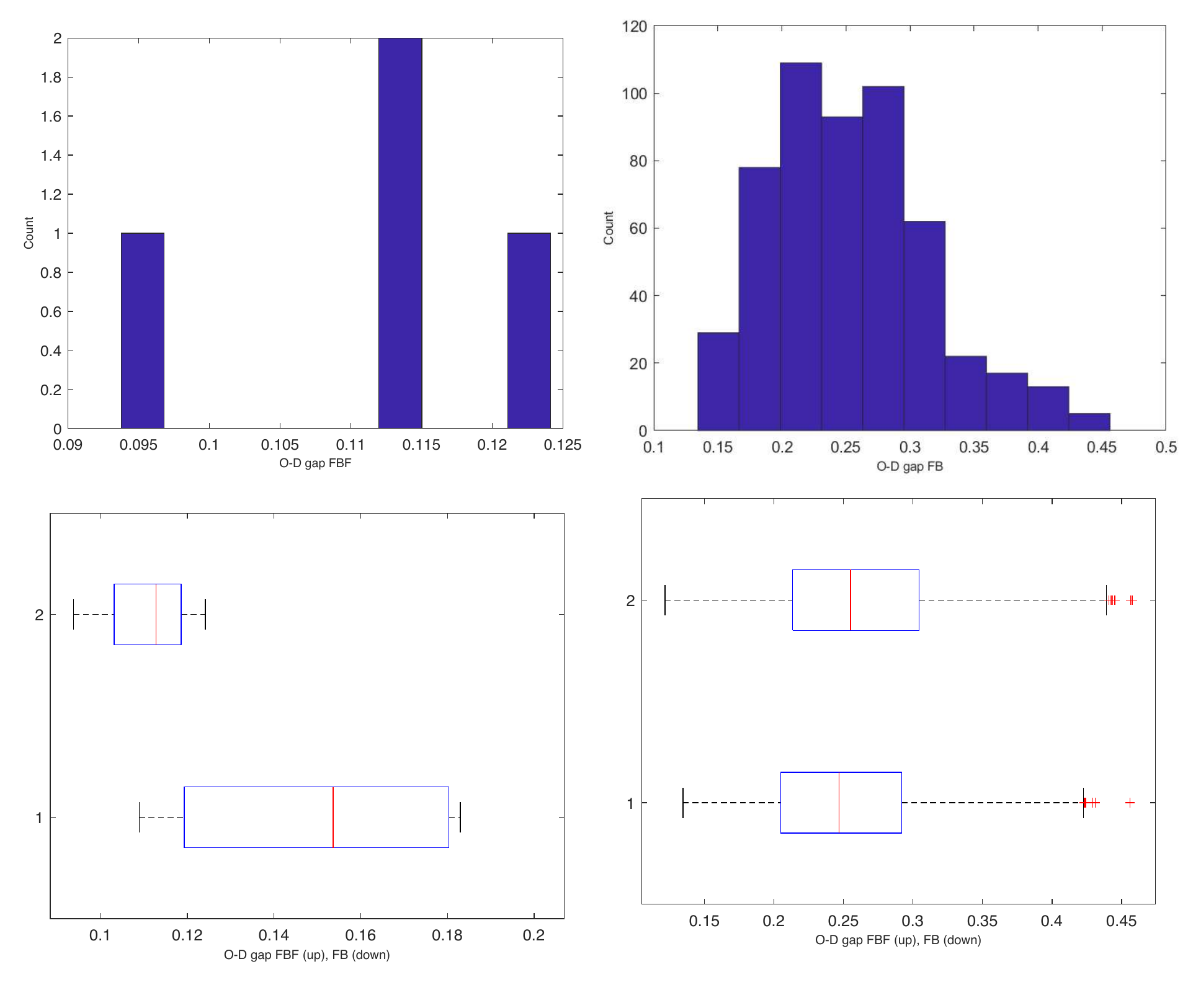}
    \caption{Distributions of O-D gaps corresponding to the DUE solutions. The O-D gap is calculated according to \eqref{eq:distance}.}
    \label{fig:ODgap}
\end{figure}

We implemented Algorithm \ref{alg:FBF-DUE} in MATLAB, building on the open-source MATLAB package described in \cite{HanEveFri19}.\footnote{This routine freely available under https://github.com/DrKeHan/DTA.} As DNL subroutine a numerical implementation of the LWR model is used, generating the delay operator $\Psi(h)$ at flow profile $h\in\scrH$.  By adapting this toolbox to Algorithm \ref{alg:FBF-DUE}, we compute dynamic user equilibria for the Nguyen and the Sioux fall network (see Figure \ref{fig:networks}) and compare our results with the projected gradient method. The parameters $\alpha_{k},\beta_{k}$ and $\gamma$ were chosen for each instance separately to guarantee the best convergence. 
The Nguyen network is a traffic network with 13 nodes connected by 19 links, and 4 o/d pairs. There are 24 paths to compute. The Sioux fall is a significantly larger instance, consisting of 76 links, 24 nodes, 530 o/d pairs and 6,180 paths. We stop the algorithm if the relative gap is smaller than a user defined tolerance, i.e.
\begin{equation}\label{eq:gap}
\eps_{k}:=\frac{\norm{h^{k+1}-h^{k}}^{2}}{\norm{h^{k}}^{2}}\leq 10^{-4}.
\end{equation}
This measure can be interpreted as the iteration complexity of the algorithm employed. Figure \ref{fig:gaps} shows the relative gaps for the Nguyen and the Sioux fall networks until the stopping criterion is reached.
It can be seen from this Figure that both methods have a similar iteration complexity, with a slight tendency favoring our FBF approach. Figure \ref{fig:paths} shows the path departure rates as well as the corresponding effective path delays. We observe that the departure rates are nonzero only when the corresponding effective delays are equal and minimum, which conforms to the notion of DUE. 
To rigorously assess the quality of obtained DUE solutions, we define the gap function between each o/d pair $w\in\scrW$ as 
\begin{align}
\nonumber
\Gamma_{w}=&\max\{\Psi_{p}(h^{\ast},t),t\in[t_{0},t_{1}],p\in\scrP_{w}\text{ such that }h^{\ast}_{p}(t)>0\}\\
&-\min\{\Psi_{p}(h^{\ast},t),t\in[t_{0},t_{1}],p\in\scrP_{w}\text{ such that }h^{\ast}_{p}(t)>0\}
\label{eq:distance}
\end{align}
In an exact DUE, we should have $\Gamma_{w}=0$ for all $w\in\scrW$. Figure \ref{fig:ODgap} displays histograms of o/d gaps obtained by running FBF and the projection method of \cite{HanEveFri19} until the stopping criterion is reached. It is seen that most o/d gaps are varying between 0.1 and 0.3 for both test instances, reflecting the early stopping of the method. We highlight that Algorithm \ref{alg:FBF-DUE} beats the projection method in the Nguyen network significantly, while it is comparable in overall performance in the Sioux fall network, and at the same time is a strongly convergent method. This provides strong evidence for the good performance of our scheme.
\section{Conclusions and Perspectives}
\label{sec:conclusion}
In this paper, we developed a new strongly convergent numerical scheme for Hilbert-space valued variational inequality problems. We implemented our algorithm in order to solve a challenging class of dynamic user equilibrium problems, and verified its competitiveness with state-of-the-art solvers used in the transportation science literature. It seems to be possible to extend our scheme to a larger class of variational problems, where distributed implementations are important, such as generalized Nash equilibrium. We leave these issues for future research. 

\section*{Acknowledgments}
This work has been completed during the thematic program LFB2019 "Maximal Monotone Operator Theory" at the ESI Vienna. The work by M. Staudigl is supported by the COST Action CA16228 "European Network for Game Theory". D. Meier and T. Vuong acknowledge financial support from the Austrian Science Fund (FWF), project I 2419-N32 and project M2499-N32 and the Doctoral Program Vienna Graduate School on Computational Optimization (VGSCO), project W1260-N35.


\end{document}